# Women mathematicians in France in the mid-twentieth century[*]

Yvette Kosmann-Schwarzbach

A short outline of the French system of "Écoles normales supérieures" and "agrégations" serves as the introduction to our account of the careers of the five women who completed a doctorate in mathematics in France before 1960 and became internationally known scientists: Marie-Louise Dubreil-Jacotin (1905-1972), Marie-Hélène Schwartz (1913-2013), Jacqueline Ferrand (1918-2014), Paulette Libermann (1919-2007) and Yvonne Choquet-Bruhat (b. 1923). This is followed by a more general description of the place of women on the mathematical scene in France between 1930 and 1960, together with a brief sketch of the accomplishments of some other women and the identification of all those who were active in research before 1960 and became professors in the French university system.

## Introduction

The intersection of the fields of history of mathematics and history of women is notoriously small. From Hypatia of Alexandria to Emmy Noether, very few women are known to have contributed to the development of mathematics, and the number of those born in France is even smaller. Were there any before Gabrielle-Émilie Le Tonnelier de Breteuil, marquise Du Châtelet-Lomond (1706-1749), the little-known Nicole-Reine Étable de Labrière Lepaute (1723-1788) and the even lesser known Marie Anne Victoire Pigeon d'Osangis (1724-1767)? Were there any between Sophie Germain (1776-1831) and those whose names and accomplishments will appear in this article? There were two French lady scientists working with the astronomer Jérôme de Lalande, Louise du Pierry (1746-1806) and Jeanne Lefrançais de Lalande (1769-1832). No other woman astronomer is known to have been active in France until Dorothée Klumke (1861-1942), who was born in the United States, had studied mathematics and astronomy at the Sorbonne, and became the third woman in 1893 to be awarded a doctorate of science in France, and the first in astronomy.

It is well known that the number and importance of women mathematicians in France in the twentieth century was due in part to the existence of the "École normale supérieure de jeunes filles", commonly called "École de Sèvres" from its former location, just outside Paris, which had a separate competitive entrance examination from that of the men's counterpart, the "École normale supérieure", whose location was and still is 45 rue d'Ulm in Paris, whence the "École normale supérieure de la rue d'Ulm". While the first École normale was founded during the Revolution in 1794 to train teachers for the schools of the new republic, the École for young women was created more than a century later, under the Third Republic, in 1881, and was not part of the national higher education system until 1936. The students took courses given by Sorbonne professors who came to teach the "sévriennes". They were finally granted the right to matriculate simultaneously at the University in 1939. Admission of women to the "rue d'Ulm" was exceptional but permitted until 1936, and still possible for two more years. Similarly, the "agrégation" competitive examination, a high-level high-school teaching certificate, was separate for men and women, but nevertheless women were permitted to take the men's examination until 1941.In 1935, it was decided to unify the requirements for the entrance examinations to the two Écoles, but for four more years, "as a transitory measure", the women's requirements remained more limited than the men's. The women studied for the entrance examination for Sèvres in a few "classes préparatoires" in high schools for girls, the best-known one for the Paris area being in the Lycée Fénelon, in the Latin Quarter. Those women who intended to sit for the competitive examination of the rue d'Ulm spent the year of preparatory studies as a tiny minority in a so-called "classe de spéciales" in a high school for boys. (These preparatory studies are at the level of the first two years of college, but they have traditionally been taught in high schools.) The two Écoles, one for men and one for women, were finally united in 1985 and it is a fact that the number of women admitted annually to the mathematics and physics programmes has decreased significantly since then. During the eight years that followed the "fusion" of the two Écoles, it happened twice that only one woman was admitted to the mathematics programme, and no more than four women were ever admitted, out of a contingent of 40 to 45 students.

Although these "écoles normales" had been created as teacher training colleges for boys' schools in the one case and for girls' schools in the other, these Écoles, admitting only the very best students, together with the École polytechnique (in the case of men), whose *raison d'être* was the training of military engineers, have been the nursery for the education of most of the French men and women who would leave a mark on the sciences. Of the five women mathematicians whose careers are sketched below, the first three had been admitted to the "rue d'Ulm" when it was still possible for young women to participate in the competition for admission, and the other two were admitted later to the École normale supérieure de jeunes filles, where the level of instruction had risen by then.

When writing about the lives and careers of scientists who pursued their studies or their careers in France during the Second World War, it is not possible to avoid references to their fate and their conduct in the trying years of the occupation. These references do not imply any judgement on my part, they are here only to serve as landmarks for future, more thorough biographies. I have also refrained from writing more than a few facts about the private lives of the mathematicians who appear in this article. All scientists have private lives and there is no need to treat those of women scientists in greater detail than those of male scientists.

# Marie-Louise Dubreil-Jacotin (1905-1972)

Marie-Louise Dubreil-Jacotin was the first woman to be named professor of mathematics in any university in France. That was in 1943, at the University of Poitiers, in the capital of the Vienne department which had been occupied, as had been all the western part of the department, by the Germans since 1940 (and would remain so until 1944), and in 1956 she obtained a professorship in Paris. Her husband since 1930 was the algebraist Paul Dubreil (1904-1994), professor at the University of Nancy during the war and then at the University of Paris from 1946 until his retirement.

Mlle Jacotin had been admitted to the École normale supérieure (rue d'Ulm) in 1926. In fact, she was ranked second in the entrance examination, but because she was a woman, her actual admission was delayed by several months. The competitive examination of aggregation was then comprised of two separate examinations, one for men and one for women. She competed with the men and succeeded, taking third place *ex-æquo* in 1929 with Claude Chevalley (1909-1984), who was to become one of the best algebraists of his generation. Early in 1930, she was invited to Oslo by the physicist, Wilhelm Bjerknes, who was working on the oscillations of warm and cold fronts in the atmosphere. She then travelled with her husband to Italy and Germany, and her meeting with Tullio Levi-Civita (1873-1941) in Rome in 1930 determined her choice of a research topic in fluid mechanics, while their meeting with Emmy Noether (1882-1935), first in Frankfurt where Noether was teaching a one-semester course in 1930-1931, then in Göttingen the following year, introduced them to "modern algebra" and re-oriented her research towards algebra. She defended a doctoral thesis, "Sur la détermination rigoureuse des ondes permanentes périodiques d'ampleur finie" (On the rigorous determination of periodic permanent waves of finite amplitude), in 1934, written under the direction of Henri Villat (1869-1972), that was soon published as an article with the same title in the *Journal de Mathématiques pures et appliquées*. Only one woman before her had defended a thesis in pure mathematics in France, Marie Charpentier, of whom more below, at the University of Poitiers in 1931, while the astronomer Edmée Chandon (1885-1944), who had successfully passed the agrégation (for women) in mathematics in 1908, had become the first woman to obtain a doctorate of mathematical sciences in France in 1930, with a thesis in astronomy and geodesy, "Recherches sur les marées de la mer Rouge et du golfe de Suez" (Researches on the tides of the Red Sea and the Suez Gulf).

Mme Dubreil-Jacotin taught the prestigious Cours Peccot at the Collège de France in 1935, a distinction awarded every year to a recent doctor in mathematics. She first taught at a university in Rennes in 1939, then in Lyons in 1940, before becoming a full professor in Poitiers, in the difficult conditions of the war years, and later in Paris. In addition to having written eight research papers and numerous notes in the *Comptes rendus de l'Académie des*

*Sciences*, first on fluid mechanics, then on algebra, she was the co-author of two books, a monograph on lattice theory and a textbook, *Leçons d'algèbre moderne*. Her research monograph, published in 1953, was the result of a collaboration with Léonce Lesieur (1914-2002), who became professor at the University of Paris in 1960, later at Orsay (Paris-Sud), and Robert Croisot (1922-1966), who became professor in Besançon. Croisot's doctorate, also published in1953, dealt with semi-modular lattices and he had written three articles with Dubreil-Jacotinin the preceding two years, while she returned to publishing, alone or in collaboration, after a lull in her scientific production. Her textbook, written with Dubreil, made Emmy Noether's abstract algebra accessible to students; it went through two editions and was translated into English and Spanish. In Paris, there was a seminar on algebra and number theory, directed since 1946 by Dubreil and Charles Pisot (1910-1984), a well-known number theorist born in Alsace who unwittingly spent the years 1940-1946 teaching in German universities, been appointed at the University of Bordeaux upon his return to France, and was professor in Paris after 1955. Dubreil-Jacotin joined her husband and Pisot as co-director of this seminar in 1957, and Lesieur joined them in 1967. For ten years (1962-1971), Dubreil-Jacotin participated in the editorial work for the seminar proceedings, which appeared as an internal publication of the "Secrétariat mathématique". Her later work was on ordered semi-groups and, after the student revolution of1968,she taught algebra applied to information processing. She wrote a chapter on women mathematicians in the influential book of François Le Lionnais, *Les Grands courants de la pensée mathématique,* first published in 1948. A street in Paris, in the new university district of the 13th arrondissement, bears her name.

## Marie-Hélène Schwartz (1913-2013)

Marie-Hélène Schwartz was the daughter of the mathematician Paul Lévy (1886-1971) and the wife of the mathematician Laurent Schwartz (1915-2002). She was admitted to the École normale supérieure (rue d'Ulm) in 1934. Her studies were interrupted by grave health problems, but Schwartz would not be deterred and they married in 1938. Then came the war and the occupation with its terrible dangers for Jews that forced the young couple into hiding. Her thesis, "Formules apparentées à celles de Gauss-Bonnet et de Nevanlinna-Ahlfors pour certaines applications d'une variété à *n* dimensions dans une autre" (Formulas related to those of Gauss-Bonnet and Nevanlinna-Ahlfors for some mappings from a manifold of dimension *n* into another), was written while she taught at the university in Rheims after the war. It was published in 1954, partly in *Acta Mathematica*, where she thanks Georges Valiron (1884-1955) for having suggested the topic of her research and André Lichnerowicz (1915-1998) for his advice, and partly in the *Bulletin de la Société mathématique de France*. Lichnerowicz, until his retirement from the Collège de France in 1986, would advise many other women, direct their doctoral theses and encourage their research.

She was appointed to the Faculté des Sciences de l'Université de Lille in 1964 where she remained until her retirement. There, in addition to teaching first-year courses in large lecture halls to young students who did not always appreciate her pedagogy, MmeSchwartz taught advanced courses and had a small group of very successful doctoral students. She wrote several book-length lecture notes for these courses, on analytic functions and sheaf theory in 1966, on differential manifolds and transversality in 1973, on singular vector bundles in 1982, and on analytic subsets of analytic manifolds in 1988.Fibered spaces had been the subject of her research for a long time. As early as 1956, she had given a three-month set of lectures in Bogotá on "Espacios fibrados". The lecture notes, in Spanish, that she wrote for this course can be found in Paris, both in the library of the Institut Henri Poincaré and in the mathematics library of the Université Pierre et Marie Curie. All the above courses were only mimeographed and remained unpublished, but between 1960 and 1992 she published a number of papers on differential geometry and analytic spaces, their singularities and their characteristic classes, that were later recognized to be fundamental, and she continued to work and to publish after her retirement in 1981. A conference in her honour was organized in Lille in 1986, and the proceedings were edited locally, under the title "Geometry Days in honour of Marie-Hélène Schwartz".  During a day of lectures organized in Paris to mark her 80th birthday, she gave a talk on her recent work that lasted almost two hours, and two years later she delivered a lecture in the very prestigious Bourbaki seminar. She wrote two research monographs that appeared in 1991 and 2000, that summarized and systematized her work on analytic stratifications and on the Chern classes of singular sets, while also confronting it with other approaches. She played a remarkably modest but efficient role as the wife of a great mathematician who conducted manifold political and social activities, so her own achievements tended to be overshadowed by his enormous mathematical talent. But her pioneering role in the field of analytic geometry is now well established.

Jacqueline Ferrand (1918-2014)

Jacqueline Ferrand, who died on 26 April 2014, was admitted to the École normale supérieure (rue d'Ulm) in 1936. Immediately after she obtained her agrégation in 1939, having been allowed to compete in the examination for men and having succeeded brilliantly, first *ex-æquo*, she became an "agrégée préparatrice"(a tutor, colloquially called a "caïmane") at the École normale supérieure de jeunes filles. In this position, she tutored the students, attempting to bring the level of Sèvres up to that of the rue d'Ulm, preparing them for the agrégation examination (for women) that they would take after three years in the École, and encouraging them in their mathematical endeavours.

In 1942, she defended her thesis entitled "Étude de la representation conforme au voisinage de la frontière" (A study of the conformal representation in the neighbourhood of the boundary), written under the guidance of Arnaud Denjoy (1884-1974), andfor which she was awarded the Girbal-Barral prize by the Académie des Sciences in 1943. Then, in 1946, she taught the Cours Peccot at the Collège de France, as had Marie-Louise Dubreil-Jacotin eleven years before. Since teaching in one or several provincial universities, before obtaining a position in Paris, was customary for all new doctors of science, she was first appointed assistant professor in 1943 at the University of Bordeaux, which was in the occupied zone until the end of August 1944,then as a professor in Caen in 1945,and in Lille, from 1948 to1956. This is when she developed the topological concept of thinness of a set at a boundary point.

She was named to a chair at the Faculté des Sciences de l'Université de Paris in 1956 – the same year as Dubreil-Jacotin, four years before Yvonne Choquet-Bruhat and ten years before Paulette Libermann –, where, in 1960, she succeeded the famous mathematician Henri Cartan (1904-2008) as the teacher of the advanced undergraduate course on functions of a complex variable, which was taught in the large "amphithéâtre Hermite" at the Institut Henri Poincaré. During her marriage to the mathematician Pierre Lelong (1912-2011), she published first under the name Jacqueline Lelong, then Lelong-Ferrand, but reverted to her maiden name after 1977.She devoted much time to the writing and publishing of several series of textbooks, corresponding to the complete undergraduate mathematics curriculum, including introductory modern differential geometry, that went through many editions and became a standard in the mathematical education of several generations of students. Meanwhile she pursued her research in complex geometry which was important enough for her to be invited to speak at the International Congress of Mathematicians in Vancouver in 1974.In 1971, she had proved that whenever a compact Riemannian manifold of dimension greater than 2 is not conformally equivalent to a sphere, its group of conformal transformations is compact, along-standing conjecture of Lichnerowicz, which Christiane Barbance (b. 1937), a former student at the École normale supérieure de jeunes filles and a doctoral student of Lichnerowicz, had refined in her thesis in 1969.

Jacqueline Ferrand retired in 1984but maintained a busy schedule of publication until 2000, and kept herself informed of developments in the field of conformal geometry, for many years afterwards.

## Paulette Libermann (1919-2007)

Paulette Libermann, one of the three daughters of a modest Parisian family with roots in Russia and Ukraine, entered the École normale supérieure de jeunes filles de Sèvres in 1938. Since 1936, this École had been headed by the physicist Eugénie Cotton, née Feytis (1881-1961), who had been a student in 1901 at the time that Marie Curie (1867-1934), Jean Perrin (1870-1942) and Paul Langevin (1872-1946) taught there, then an "agrégée répétitrice"(a position later called "agrégée préparatrice") at the École from 1905 to 1931. She had defended a thesis on magnetic salts in 1925, would later work on the magnetic properties of rubber, and she is remembered for her very active social and political role in feminist and communist circles. Mme Cotton was responsible for the reform that brought the École de Sèvres into the higher education system in 1936 and, while the reform itself took several years to be fully implemented, as early as 1939 she hired the gifted and enterprising young mathematician Jacqueline Ferrand and entrusted her with the tutoring of the students. New teachers were invited to lecture at the École. Several were "young and brilliant like André Lichnerowicz, or world renowned like Élie Cartan", as Libermann wrote in [1995, p. 19]. Lichnerowicz had been appointed in 1938 and Élie Cartan (1869-1951) would teach at Sèvres after his retirement from the Faculty of Science in Paris in 1940.

The war broke out in September 1939, the German army invaded a large part of France that included Paris, and there followed in the spring and summer of 1940 a southward exodus of civilians. Paulette and her family took refuge in Lyons, where she was able to sit for two examinations. She then travelled to Toulouse to sit for additional examinations. She returned to Paris with her family in autumn. Part of the École had moved to Paris after the German troops occupied the building in Sèvres, and she prepared for the agrégation while living with her family. But she was prevented from sitting for this examination, that would have led to an appointment as a teacher in a state high school, by the "statut des Juifs" decreed by the Vichy regime in October 1940, and reinforced in 1941. Thanks to Eugénie Cotton, the courageous head of the École, she was nevertheless permitted to start doing research during the academic year 1941-1942 under the guidance of Cartan, the foremost geometer in France. He had taught the "sévriennes" at Sèvres and he now lectured for them in Paris. She wrote an essay, "Sur les propriétés projectives des courbes gauches et des complexes de droites. Leur interprétation dans l'espace à cinq dimensions et dans le plan" (On the projective properties of skew curves and line complexes. Their interpretation in five-dimensional space and in the plane), for the "diplôme d'études supérieures de mathématiques". (Cartan's report on her essay is conserved in the Fonds Élie Cartan in the Archives of the Académie des Sciences in Paris.) After the Jews were obliged to wear the yellow star and thousands were deported in the summer of 1942, Paulette and her family fled Paris a second time for Lyons, where they narrowly escaped deportation. After Lyons was liberated by the Allied forces in September 1944, she was able to return to Paris which had just been liberated. She was finally allowed to sit for the examination for the aggregation (for women) at the end of 1944 and took second place. She first taught for a trimester in 1945 in the "Lycée de Jeunes Filles" (high school for girls) in Douai in northern France, but scholarships of the British Council, and later of the Centre National de la Recherche Scientifique (CNRS) enabled her to pursue her research from 1945 to 1947 under the guidance of J.H.C. Whitehead (1904-1960) in Oxford, where she also obtained the title of bachelor of science.

She continued her research in geometry, this time in Strasbourg. She first taught in the high school for girls from 1947 to 1951. Then, for the next three years in Strasbourg she held a research position funded by the Centre National de la Recherche Scientifique which enabled her to complete her doctoral thesis, "Sur le problème d'équivalence de certaines structures infinitésimales" (On the equivalence problem for some infinitesimal structures). She defended it in 1953, having already co-authored several papers with her thesis advisor, the distinguished geometer Charles Ehresmann (1905-1979). It appeared in the same year in *Annali di Matematica*.

She was to be the first "sévrienne" to become an important mathematician and a full professor of mathematics in any university. Her first position in a French university was as an assistant professor at the University of Rennes, in 1954, where she was promoted to the rank of full professor in 1958. She spent a semester as an invited professor at the University of São Paolo in 1961.Then, in 1966, she was named to a chair in Paris. When the Faculty of Science of the University of Paris split after 1968, mostly along ideological lines, she chose Paris VII in accordance with her left-leaning convictions. A small scientific conference marked her retirement in 1986, and a week-long international conference was organized in Paris in 2009 to honor her memory and mathematical legacy.

Libermann is the author with Charles-Michel Marle (b. 1934), professor at the University of Paris VI, of an important text on symplectic geometry and its applications to mechanics, first printed as four volumes of the *Publications de l'Université Paris VII* in 1986 and published in English translation in 1987, which has become a classic. She published numerous research papers on many different aspects of differential geometry and global analysis (contact manifolds, foliations, connections, groupoids and Lie algebroids, to name just a few). She was invited to prestigious universities and conferences, she travelled extensively, and was active in research well after her retirement and until a few months before her death. In 1996, she wrote a short scientific biography of Élie Cartan, and remained devoted to the Cartan family all her life. The collection of her mathematical papers and books is conserved in the mathematical library and archives of the Université Pierre et Marie Curie in Paris. The importance of her work, especially her publications of the 1950s, is now being increasingly recognized.

# Yvonne Choquet-Bruhat (b.1923)

Yvonne Choquet-Bruhat is the daughter of the physicist Georges Bruhat (1887-1945), acting director of the École normale supérieure in 1941, who was deported to Buchenwald by the Nazis for resisting the Gestapo and died in the Sachsenhausen camp. Her brother was the algebraist François Bruhat (1929-2007), professor at the University of Paris from 1961 until his retirement in 1989. She was a student at the École normale supérieure de jeunes filles, 1943-1946, took first place in the aggregation examination in 1946, and was soon appointed as "agrégée préparatrice" at the École. In her position as a research associate at the Centre National de la Recherche Scientifique, she completed her doctoral thesis in 1951, on the hyperbolic systems that need to be solved in the general theory of relativity, "Théorèmes d'existence pour certains systems d'équations aux derives partielles non linéaires" (Existence theorems for some systems of nonlinear partial differential equations). It appeared in *Acta Mathematica* the following year. Her advisers had been Georges Darmois (1888-1960), Lichnerowicz and Jean Leray (1906-1998). Darmois had published on general relativity from 1923 to the early 1930s, although his later work was mostly in probability theory and statistics. Lichnerowicz, who had written his thesis on the Einstein equations in 1939 as Darmois' doctoral student, suggested that she should generalise the 3+1 decomposition that Darmois and himself had obtained in particular coordinates. Leray suggested that she rather study the non-analytic Cauchy problem for the Einstein equations, a topic on which she obtained fundamental results for evolution, related to the existence of gravitational waves, using a method initiated by Sergei Sobolev (1908-1989). She obtained results on constraints for initial data in subsequent publications. Leray, who in the early 1950s was working on general hyperbolic problems, oriented her research towards the analysis of partial differential equations on manifolds. At the time, and for some years afterwards, he spent every fall semester at the Institute for Advanced Study in Princeton. He invited her to visit the Institute, and she arrived in 1951, shortly after the defense of her thesis. In a recent interview, she recalled that Leray had introduced her to Albert Einstein, telling him that she had worked on general relativity and that her father was the physicist Bruhat who had died in deportation. Einstein asked her to give him an account of her thesis, she did, and he found it interesting.

She left Princeton for France in 1953 to take up her first university position, associate professor at the University of Marseilles, where she remained until 1958, with the exception of another visit to Princeton in 1955.Einstein died in April of that year and she always regretted not having visited his office more often during her earlier stay in Princeton, since he had invited her to stop by as often as she would like. After Marseilles, she taught at the University of Rheims during the year1958-1959, after which she was elected professor at the Faculté des Sciences of the University of Paris in 1960, where she remained until the year of her retirement, 1992.

She has published continuously since 1948 and is the author of nearly 250 scientific papers dealing with analysis, differential geometry, general relativity and mathematical physics, and of several books, in particular and most recently, an 800-page volume on the solutions of the Einstein equations. She was the first woman to be elected to the French Académie des Sciences, first as a correspondent in 1978, and then as a member in 1979. (The first woman to become a correspondent was the physicist Marguerite Perey (1909-1975). Choquet-Bruhat was the first woman in any field of science to become a member. The second woman mathematician was MichèleVergne (b. 1943) in 1997, and the third was Claire Voisin (b.1962)in 2010.) Her husband was the mathematician and professor at the University of Paris, Gustave Choquet (1915-2006). She also published under the name Fourès-Bruhat, and under her maiden name, Bruhat.

Choquet-Bruhat recalls that her father taught her to love physics, but that he did not envision a scientific career for her, rather the life of a good mother and high school teacher. She became an internationally recognized mathematician, whose work has important applications in astrophysics, in particular in the study of gravitational waves. She has been awarded many honours, both French and American, medals and prizes, and in 2008 she was named "grand officier" in the Légion d'Honneur. An international colloquium in her honour, "Physics on Manifolds", was held in Paris in June 1992, when she became professor emeritus, and a collection of articles with the same title was published as a tribute to her contributions to both physics and geometry. On the occasion of her 90th birthday in January 2014, a prestigious international gathering at the Institut des Hautes Études Scientifiques, in Bures-sur-Yvette near Paris, where she enjoyed having a desk after her retirement, celebrated her work.

Personal endnote

While I was a student in the early 1960s, three of the professors teaching the courses in mathematics for the undergraduate and beginning graduate students at the Institut Henri Poincaré of the University of Paris or at the newly completed buildings of the Faculté des Sciences along the quai Saint-Bernard were women, which was remarkable for the time. Two were my teachers, Mme Lelong, as we called her, who taught functions of a complex variable, and Mme Choquet-Bruhat who taught an advanced course on general relativity, while Mme Dubreil was teaching an algebra course which I did not take. That was the time when Louis de Broglie (1892-1987) was delivering his lectures on wave mechanics in the "Amphithéâtre Darboux" of the same Institut, for the last time, because he retired the following year, 1962. Simultaneously, in the laboratoire de physique théorique on the second floor of the Institut, there were other women scientists working on relativity theory, Mme Tonnelat (Marie-Antoinette Tonnelat, 1912-1980), who was rumored to have been an assistant to Einstein, had taught the Cours Peccot in 1944-1945, and had worked with Schrödinger in Dublin, *circa* 1946, and her own assistant, Mlle Mavridès (Stamatia Mavridès, b. 1924), who pursued her career in cosmology at the Meudon branch of the Paris Observatory.

In 1966, when I was a beginning researcher, Paulette Libermann arrived in Paris. We met at the "thé des mathématiciens", a social gathering that took place every Wednesday afternoon in a large room on the first floor of the Institut Henri Poincaré, and she befriended me. (I think that she saw me as a younger version of herself, an unmarried female mathematician.) She soon invited me to participate in her seminar and to give a talk, which gave me my first opportunity to speak about my work. (This seminar was attended by a small number of participants, yet my talk took place in the impressive "Amphithéâtre Hermite" in which I had been a student a few years before.) Much later, from 1982 to 1990, we organized a seminar together. Then we became friends.

When I was first appointed professor at the University of Lille in 1969, I became a colleague of Marie-Hélène Schwartz, and she welcomed me warmly. My office was next door to hers and, later, when I told her that I was engaged to get married, she gave me precious, maternal advice. (Maybe she saw me as a mathematical daughter.) We sometimes held animated discussions on the train that would take us back to Paris where we both lived, after a day of teaching in Lille. We remained friends after she retired, although we had fewer occasions to meet.

With Jacqueline Ferrand, several decades after I had been her student, I collaborated when it was, sadly, time to write an obituary for Paulette Libermann.

After participating in many of the same conferences, I can now consider my former teacher Yvonne Choquet-Bruhat, whose work I much admire, to be a friend.

I shall write a few words about four other French women in science whom I have known personally.

Lucienne Félix (1901-1994) was a student at Sèvres from 1920 to 1923, and, after teaching in a high school in Lille, returned as an "agrégée répétitrice", then also called "professeur adjointe de mathématiques", from 1929 until 1938. When the war broke out, she was teaching at the college level in the girls' high school in Versailles, preparing them for the entrance examination for Sèvres, but she was dismissed from her teaching position in 1940 because she was Jewish. Henri Lebesgue (1875-1941) had been her teacher at Sèvres and she had served as an assistant to him and to Henri Villat, when she taught there. Lebesgue, who was already ailing when he taught his course at the Collège de France in 1941, invited her to audit his lectures in spite of the restrictions that applied to Jews, and he entrusted his preparatory notes, a copybook and notebooks, often handwritten in pencil, to her. For three years, she organized these elements into a book, even writing up a chapter that Lebesgue had never lectured upon. In the spring of 1944, she was able to return Lebesgue's manuscripts to his family and to deliver her completed manuscript to Paul Montel (1876-1975), a member of the Academy of Sciences, who disdainfully "asked [her] to leave it with the concierge" [Félix 2005, p. 66]. When she was arrested with her father a few months later by German policemen, "[her] first thought was that [she] had fortunately put Lebesgue's work out of harm's way"[pp. 66-67]. Interned in Drancy, she and her father escaped being deported in the last trains to the concentration camps. Lebesgue's "Leçons sur les constructions géométriques" (Lectures on geometric constructions) were published posthumously in 1950, with a preface by Montel that does not acknowledge her role in the completion of the book's manuscript. He obviously had never read it, since he did not correct an enormous, glaring misprint, which she immediately saw upon receiving the printer's proofs in 1949. Lucienne Félix was reinstated in a high-level teaching position after the war, only to be again fired by an anti-semitic principal. After various temporary teaching jobs, including teaching a "classe de [mathématiques] spéciales" for young men who had fought in the war or in the resistance, she was finally appointed to a Parisian high school for girls where she made a point of teaching at every level. She was very active in research in the pedagogy of mathematics, was a pioneer in the introduction of "modern mathematics" at the high school level, and even in elementary schools, about which she published many books and articles, and she lectured all over the world. She is also remembered as an active supporter of the Russian mathematicians who were *refuznik*s. (MlleFélix was my teacher when I was 12 and again in my last year of high school. At that time, not only did she encourage me, but she met my father and insisted that I ought to go on and study in a "classe préparatoire".)

Cécile DeWitt-Morette (b.1922), who studied during the war at the University of Caen and then at the Sorbonne in Paris, ended her career as both a member of the governing council of the Institut des Hautes Études Scientifiques and a professor with an endowed chair in the physics department of the University of Texas, Austin, where her husband, the physicist Bryce DeWitt (1923-2004), was a professor. When she was a student in 1944, in spite of the rumors of the imminent landing of the Allied troops, she courageously left Caen for Paris to sit for an examination scheduled for June 6, and this saved her life since, on that day, her entire family was killed in the bombing of their home in Caen. Much later, in Texas, rules that prohibited the hiring of husband and wife in the same university prevented her for 16 years from obtaining a suitable position, but not from publishing and directing numerous graduate and Ph. D. theses. She is the author of over a hundred papers in mathematical physics and is the co-author, with Yvonne Choquet-Bruhat, of a two-volume course on analysis and geometry which went through seven editions. Most recently (2006), she wrotea book on functional integration with the mathematician Pierre Cartier (b. 1932), in which varied advanced mathematical methods are used to treat, among other topics, path integration, supersymmetry, stochastic processes and quantum field theory.

I must also evoke the memory of my friend Claudette Buttin (1935-1972), née Vanmaele, who had entered the École normale supérieure de jeunes filles in 1954, worked in differential geometry under the direction of Lichnerowicz, but did not live to defend her thesis. An article in the *Bulletin de la Société mathématique de France* which describes her work was drawn up from her notes by Pierre Molino (b. 1935) and published in 1974 under her name.

Yvette Amice (1936-1993) entered the École normale supérieure de jeunes filles in 1956, obtained the agrégation in 1959, defended a thesis on *p*-adic interpolation in 1964, and taught the Cours Peccot in 1966. That year she was named professor at the University of Bordeaux and in 1970 she was appointed to the Université Paris VII where she directed an active seminar in number theory until her premature death from cancer. She had many varied activities in research and administration. She was one of three vice-presidents of the Société mathématique de France in 1974 and the only female member of its 25 member council, and in 1975 she became its second woman president. (Dubreil-Jacotin had been president in 1952.) It was not until 1998 that another woman mathematician, the algebraic geometer Mireille Martin-Deschamps, who had entered the École normale supérieure de jeunes filles in 1965, presided over the national professional organization of mathematicians, this time for three years, the statutes having changed. Two other women have been president for three years since then.

Who else?

In ÉlieCartan's family, his son Henri became internationally known as a topologist, and two of the women were also gifted in mathematics, although their achievements do not bear comparison to those of these two giants of twentieth century mathematics.

Anna Cartan (1878-1923) was Élie's younger sister. She entered the science section of the École normale supérieure de jeunes filles de Sèvres in 1901, together with only three other young women, one of them Eugénie Feytis. She then studied mathematics and succeeded in the agrégation for women in 1904. Lucienne Félix recalls the efficient and perceptive role of Mlle Cartan as an instructor of training sessions for future teachers of mathematics at the École, that were offered to the "sévriennes" after they obtained the aggregation but before they went out to one provincial town or the other to teach in high schools.

Hélène Cartan (1917-1952) was the youngest of Élie Cartan's four children. It was to the École normale Supérieure (rue d'Ulm) that she was admitted in 1937, before the reform barred women from sitting for the entrance examination to the men's École. In 1940, as she was ready to sit for the agrégation, only the examination for women was organized because of the war. She sat for it and took first place. She became a high school teacher, but succeeded in publishing a note in the *Comptes rendus de l'Académie des Sciences* in 1942. Later however, she became ill with tuberculosis and was unable to do research or even to teach.

The physicist Jeanne Lattès, née Ferrier (1888-1979) had studied mathematics and physics in Montpellier, and she defended a thesis in physics in 1926. In 1930, because of her declining health, she had to leave her position at the Curie Laboratory. She then moved to the nearby Institut Henri Poincaré where she assisted Émile Borel (1871-1956) in his work on probability theory.

When Paulette Libermann was appointed to the University of Rennes in 1954, there was already one woman on the faculty of the department of mathematics, Marie Charpentier (1903-1994), and there would be a third four years later. Mlle Charpentier had joined the Société mathématique de France in 1930, an indication that she was then considered an active researcher. She may have been only the second such female member. (Edmée Chandon had become a member in 1919, Lucienne Félix would become a member in 1946.) The following year, Marie Charpentier defended her thesis in mathematics in Poitiers before a jury chaired by the Paris professor, Montel, "Sur les points de Peano d'une equation différentielle du premier ordre" (On the Peano points of a first-order differential equation). She was a Rockefeller fellow at Harvard University in 1931-1932 where she studied and worked with George Birkhoff (1884-1944). One assessment of her personality and work quoted by Reinhard Siegmund-Schultze [2011, p. 125] affirmed that "she has great courage and is up to the average of men fellows in ability". But still, when she returned to France, she had to take a teaching position in a high school for lack of any available research position. It was only when the government of the Front Populaire came to power in 1936 and the prime minister, Léon Blum (1872-1950), had chosen the Nobel prize winner Irène Joliot-Curie (1897-1956) as under-secretary for research, one of three women under-secretaries in his government, that Charpentier received a research grant from the Caisse Nationale de la Recherche Scientifique, which had been created the year before. She published some 15 articles and notes in the *Comptes rendus de l'Académie des Sciences* between 1930 and 1939, and she also sat for the agrégation competition for women in 1936 and took third place. She was finally appointed to the University of Rennes in 1943, after which she published only one more paper, in 1946, on the topological properties of solutions of differential equations. By the time MlleLibermann arrived in Rennes, Mlle Charpentier had been promoted to the rank of professor, and she taught there until her retirement in 1973.

Huguette Delavault (1924-2003) was a former student of the École normale supérieure de Fontenay-aux-Roses to which she had been admitted in 1946. This École for women was first created in 1880 to train teachers who would themselves be responsible for the training of elementary school teachers, in parallel with the École normale supérieure de Saint-Cloud for men, created in 1882. These two Écoles evolved after 1945 to become counterparts to the Écoles normales supérieures, but on a slightly lower level. Delavault was certainly the first "fontenaisienne" to become a doctor in mathematics and finally a professor. She successfully sat for the agrégation examination in 1952, then held a research position at the Centre National de la Recherche Scientifique while preparing her thesis under the direction of Henri Villat, "Application de la transformation de Laplace et de la transformation de Hankel à la determination de solutions de l'équation de la chaleur et des équations de Maxwell en coordonnées cylindriques" (Application of the Laplace transform and of the Hankel transform to the determination of solutions of the heat equation and the Maxwell equations in cylindrical coordinates), which she defended in 1957. She had joined the Société mathématique de France in 1953 and had announced the results of her thesis in notes in the *Comptes rendus de l'Académie des Sciences*. She started teaching at the University of Rennes in 1958, published a memoir on the applications of integral transforms in several variables, and was named a professor in 1962. After that date, she ceased publishing but later held important administrative positions, in particular as assistant head of the École at Fontenay-aux-Roses from 1976 to 1980, and worked ceaselessly to further equal opportunities for women in science. She was the driving force in the creation of the association "Femmes & Sciences" in 2000.

Lille was another university where there was a relatively large contingent of women on the faculty of the mathematics department before 1960 (and many more later). They were nearly all in junior positions, but Christiane Chamfy (later Blanchard), who had entered the École normale supérieure de jeunes filles in 1951, and was the first "sévrienne" (but not the last) to be appointed in Lille, joined the Department of Mathematics in 1958, earned a doctorate at the University of Paris in 1959 and became the first woman to chair the department in 1961, before leaving for the University of Marseilles in 1964. It was only in1996 that the Department of Mathematics in Lille was again headed by a woman, Anne Duval (née Scherpereel), also a former student of the École normale supérieure de jeunes filles which she had entered in 1964.

Sévriennes and others

I have attempted to determine which other women in France published in mathematical periodicals before 1960 and became professors. Among them eight were former students at the École that was now often called by the acronym ENSJF. Below, they are listed in the order of their admission to the École, and the year of their admission follows their name. The date of their doctorate (doct.) is given within the parentheses. All these defenses were at the University of Paris, except that of Aimée Baillette which took place at the recently created University of Orsay (Paris-Sud). I have counted six women who did not attend the ENSJF but held a doctorate before 1960 and became professors. Admission to the Société mathématique de France (abbreviated to SMF below) had to be sponsored by at least one member, and was normally a sign of some research activity. I include the date of first membership.

Aimée Baillette, ENSJF 1942 (doct. 1965, professor in Perpignan), joined the SMF in 1961,
Marianne Teissier, ENSJF 1946 (later Guillemot, doct. 1969, professor in Tours), joined the SMF in 1951,

Françoise Tison, ENSJF 1947 (later Pécaut, professor in Avignon), joined the SMF in 1953,

Josette Renaudie, ENSJF 1948 (later Charles, doct. 1957, professor in Montpellier), joined the SMF in 1953,

Françoise Guyon, ENSJF 1949 (later Hennequin, doct. 1958, professor in Clermont-Ferrand), joined the SMF in 1954,

Monique Augé, ENSJF 1952 (later Lafon, associate professor at Paris XIII), joined the SMF in 1956,

Françoise Besson, ENSJF 1952 (later Bertrandias, doct. 1965, professor in Grenoble), joined the SMF in 1958,

Marthe Hugot, ENSJF 1952 (later Grandet, doct. 1964, professor in Caen), joined the SMF in 1957.

Elizabeth Lutz (doct. 1951, professor in Grenoble), joined the SMF in 1954,

Édith Mourier (doct. 1952, professor of probability at Paris VI), joined the SMF in 1953,

Simone Marquet (doct. 1955, professor in Lille), joined the SMF in 1954,

Simone Lemoine (later Dolbeault, doct. 1956, professor in Poitiers), joined the SMF in 1953,

Denise Huet (doct. 1959, professor in Nancy), joined the SMF in 1958.

Andrée Bastiani (later Ehresmann), who first published in 1958, became professor in Amiens.

Marie-Paule Brameret (later Malliavin) began publishing in 1960, joined the SMF in the same year, earned a doctorate in 1965 and became professor at Paris VI.

I estimate that less than 50 women were members of the SMF before 1960. Of those that I did not already mention, and who were living in 1974, three were retired teachers (one of them, Odette Mongeaud-Devisme, had joined in 1932, just after Marie Charpentier, Geneviève Guitel joined in 1949 and Germaine Ayrault in 1952), two had technical positions at the CNRS (Denise Lardeux and France Fages), three taught at a university (Marie-Thérèse Estival in Besançon, Gisèle Vors at Orsay, and Françoise Benzécri-Le Roy in Rennes), while Judith Winogradzki was professor of theoretical physics in Rouen and Rose Bonnet had been an astronomer at the Paris Observatory. The number of members of the Society increased rapidly after 1960, and women formed a little over 10% of the 1166 members in 1974.

The number and contents of the mathematical publications of those women who were professors varied greatly, from cases where we find only a few notes in the *Comptes rendus de l'Académie des Sciences* written when they were preparing their dissertations, and little or nothing else, to cases of a full career of publishing with over seventy publications. The fact that the French system gave tenure to teachers at the university at a very early stage in their career was a factor that was very beneficial to women, in view of their family obligations, but also allowed some men and women to continue teaching even if they were producing practically no research. It would seem that, for some of the women who married, their duties as mothers prevented them from achieving more in the way in research. Nevertheless, several women did succeed in combining marriage, maternity and mathematics.

All preparatory courses and all establishments of higher education in France have now become open to women. Among the mathematicians born during and after the Second World War who started publishing in the 1960s and later, there have been many talented women whose contributions have brought further distinction to the French scientific community, but their careers belong to a later period in the development of the role of women in science.

References

It is impossible, in a few pages, to give an account of the scientific work of the French women mathematicians I have cited, and I have not endeavoured to do so, nor shall I cite their numerous scientific articles that have appeared in periodicals. I only list some of the books that they have written, texts as well as research monographs.

I Books written by women mathematicians mentioned in this article

III Furtherreading